\documentstyle[11pt]{article}
\def\beq{\begin{equation}}
\def\eeq{\end{equation}}
\def\bea{\begin{eqnarray}}
\def\eea{\end{eqnarray}}
\def\beas{\begin{eqnarray*}}
\def\eeas{\end{eqnarray*}}
\def\nn{\nonumber}
\def\hy{\hbox{--}}

\begin{document}
\begin{center}
{\Large \bf 
Transformation and summation formulas for\\[3mm]
Kamp\'e de F\'eriet series $F_{1:1}^{0:3}(1,1)$
}\\[3cm]
{\bf S.N.~Pitre and J.\ Van der Jeugt\footnote{Senior Research
Associate of N.F.W.O. (National Fund for Scientific Research of
Belgium).} }\\[1cm]
Department of Applied Mathematics and Computer Science,\\
University of Gent,\\
Krijgslaan 281-S9, \\
B-9000 Gent, Belgium.\\[2mm]
E-mail : Pitre.Sangita@rug.ac.be and Joris.VanderJeugt@rug.ac.be.
\end{center}

\begin{abstract}
The double hypergeometric Kamp\'e de F\'eriet series
$F^{0:3}_{1:1}(1,1)$ depends upon 9 complex parameters. We
present three cases with 2 relations between those 9 parameters,
and show that under these circumstances $F^{0:3}_{1:1}(1,1)$ can
be written as a ${}_4F_3(1)$ series. Some limiting cases of
these transformation formulas give rise to new summation results
for special $F^{0:3}_{1:1}(1,1)$'s. The actual transformation
results arose out of the study of 9-$j$ coefficients.
\end{abstract}

\section{Introduction}

The results given in this paper are basically mathematical but the
subject area where they have arisen is physics. Therefore we
shall devote this introduction to giving some relevant
references and to describing the field in which these results
have appeared naturally. 

We shall present a number of hypergeometric series
transformation and summation formulas which were discovered when
systematically studying the 9-$j$ coefficient. These
coefficients, depending upon 9 integer or half-integer
parameters and usually denoted by
\beq
\left\{ \begin{array}{ccc} j_1 & j_2 & j_{12}\\ j_3 & j_4 &
j_{34}\\ j_{13}& j_{24} & J \end{array} \right\}
\label{ninej}
\eeq
appear in the quantum mechanical treatment of angular momentum
in physics. The 9-$j$ coefficient can be considered as a
transformation coefficient connecting two different ways in
which four angular momenta $j_1,j_2,j_3,j_4$ can be 
coupled~\cite{Wigner,Biedenharn}~: either $j_1$ and $j_2$ to
$j_{12}$, $j_3$ and $j_4$ to $j_{34}$, and $j_{12}$ with
$j_{34}$ to $J$; or else $j_1$ and $j_3$ to $j_{13}$, $j_2$ and
$j_4$ to $j_{24}$, and $j_{13}$ with $j_{24}$ to $J$.
By identifying the representation theory of angular
momentum with the representation theory of the Lie group
SU(2), one can also interpret the 9-$j$ coefficient as the
transformation coefficient relating irreducible constituents of
the SU(2) tensor products $(V_1\otimes V_2)\otimes (V_3\otimes
V_4)$ and $(V_1\otimes V_3)\otimes (V_2\otimes V_4)$.

From equations expressing the 9-$j$ coefficient in terms of
3-$j$ (or Wigner) coefficients (cfr.~\cite{Edmonds}),
it follows that the 9-$j$ coefficient has 72
symmetries~\cite{Wigner,Biedenharn}~: up to a sign,
the array~(\ref{ninej}) is invariant for row and column
permutations and for transposition.
On the other hand, there
exists an alternative expression of the 9-$j$ coefficient~: this
is the triple sum series expression obtained by
Ali\v sauskas and Jucys~\cite{Alisauskas} and independently
proved in~\cite{Jucys}. This triple sum series is further
discussed in Ref.~\cite{Rao}. In Ref.~\cite{Vanderjeugt}, a certain
``doubly stretched'' 9-$j$ coefficient for which a single term
expression exists was considered and it was shown how the
symmetries which do not yield a single term can give rise to
single, double and triple hypergeometric summation formulas.

Presently, we have investigated those symmetries of doubly
stretched 9-$j$ coefficients which give rise to double sums (the
reason for considering doubly stretched coefficients is that
singly stretched or unstretched coefficients would be too
general to lead to summation or transformation results). The
basic formulas for doubly stretched 9-$j$ coefficients can be
found in Refs.~\cite{Sharp} and~\cite{Varshalovich}. The clue
for finding new transformation or summation formulas lies in
comparing the expressions to which the triple sum series reduces
for two different symmetries. As such, we have systematically 
equated double sums with single sums. It turned out that all
double sums are one of the three types of Kamp\'e de F\'eriet
series $F^{0:3}_{1:1}(1,1)$ given in the next
section. In the first instance, the 9-$j$ coefficients yield only
finite series; but apparently most of the transformation or
summation formulas thus obtained have a wider region of
validity. In the following sections we shall present and prove
the three types of special Kamp\'e de F\'eriet series
transformation formulas without further reference to the
background of 9-$j$ coefficients.

\section{Transformation formulas}

The Kamp\'e de F\'eriet function $F_{1:1}^{0:3}$, a
generalization of Appell series~\cite{Appell}, is defined as
follows~\cite{Srivastava}~: 
\beq
F(x,y)=F_{1:1}^{0:3}\left[ {\hy \atop d}:{a,b,c\atop
e};{a',b',c'\atop e'};x,y\right]=
\sum_{m,n=0}^\infty {(a,b,c)_m (a',b',c')_n \over
(d)_{m+n}(e)_m(e')_n} {x^m y^n \over m!n!},
\label{kdf}
\eeq
where the notation is as in Ref.~\cite{Slater}~: $(a)_m$ is a
Pochammer symbol, and $(a,b,c)_m=(a)_m(b)_m(c)_m$.
The region of convergence is given by $|x|<1$ and $|y|<1$, and
the series is absolutely convergent for $|x|=|y|=1$
provided~\cite{Hai} 
\beq
\Re (d+e-a-b-c)>0 \qquad\hbox{and}\qquad\Re(d+e'-a'-b'-c')>0.
\label{converg}
\eeq
As usual, it is understood that no denominator parameters are
zero or negative integers.

The series $F(1,1)$ depends upon 9 parameters $a,b,\ldots,e'$.
We present three transformation formulas with 7 free parameters,
i.e.\ when 2 relations hold among the 9 parameters. The first of
our formulas is when $a'=d-a$ and $e'=d+e-a-b-c$; then we find~:
\bea
F(1,1)&=&\Gamma\left[ {d,e',d+e'-a'-b'-c' \atop
a',d+e'-a'-b',d+e'-a'-c'} \right]\nn\\[1mm]
&&\times {\ }_4F_3\left[{a,e-b,e-c,d+e'
-a'-b'-c' \atop e,d+e'-a'-b', d+e'-a'-c'};1\right].
\label{kdf1}
\eea
Herein, $\Gamma$ is the classical Gamma function with the
convention 
\[
\Gamma\left[{a_1,a_2,\ldots \atop b_1,b_2,\ldots}\right] =
{\Gamma(a_1)\Gamma(a_2)\ldots \over
\Gamma(b_1)\Gamma(b_2)\ldots}; 
\]
and ${}_4F_3$ is a generalized hypergeometric
function~\cite{Bailey,Slater}. Eq.~(\ref{kdf1}) holds
when~(\ref{converg}) and $\Re(a')>0$ are satisfied.

The second formula is for $a'=d-a$ and $d+e-a-b-c=1$.
Now we find~:
\bea
F(1,1)&=&\Gamma\left[ {e,e-a-b,e-a-c,e-b-c,e',e'-b'-c' \atop
e-a,e-b,e-c,e-a-b-c,e'-b',e'-c'} \right]\nn\\[1mm]
&&\times {\ }_4F_3\left[{a,b',c',d-b-c \atop 
d-b,d-c,1+b'+c'-e'};1\right].
\label{kdf2}
\eea
This formula is shown to be valid under one of the following
conditions~: 
\begin{itemize}
\item[(i)] $a$ is a negative integer and the second inequality
of~(\ref{converg}) is satisfied; 
\item[(ii)] $c$ and $c'$ are negative integers.
\end{itemize}

The third and last of our transformation formulas is for
$a'=d-a$ and $b'=d-b$; moreover we assume that $a'$ or $b'$ is a
negative integer. Then~:
\bea
F(1,1)&=&\Gamma\left[ {d,e,a+b-d,d+e-a-b-c \atop
a,b,e-c,d+e-a-b} \right] \nn \\[1mm]
 &&\times{\ }_4F_3\left[{a',b',e'-c',d+e
-a-b-c \atop e',1+d-a-b, d+e-a-b};1\right].
\label{kdf3}
\eea
In this case, the termination of the ${}_4F_3$ series implies
that conditions~(\ref{converg}) are sufficient. 

From the three transformation formulas (\ref{kdf1})-(\ref{kdf3})
we believe that the first two are new. The third is not new; in
fact it is the terminating form of a formula given by
Karlsson~\cite{Karlsson2} in his proof of one of our earlier
double hypergeometric summation formulas~\cite{Vanderjeugt} (see
eq.~(\ref{res1}) of this paper). Karlsson's formula holds for
$a'=d-a$, $b'=d-b$, without the assumption that $a'$ or $b'$
should be negative integral. It reads~:
\bea
F(1,1)&=&\Gamma\left[ {d,e',a'+b'-d,d+e'-a'-b'-c' \atop
a',b',e'-c',d+e'-a'-b'} \right]\nn\\[1mm]
&&\times {\ }_4F_3\left[{a,b,e-c,d+e'
-a'-b'-c' \atop e,1+d-a'-b', d+e'-a'-b'};1\right] \nn\\[1mm]
&+&\Gamma\left[ {d,e,a+b-d,d+e-a-b-c \atop
a,b,e-c,d+e-a-b} \right]\nn\\[1mm]
&&\times {\ }_4F_3\left[{a',b',e'-c',d+e
-a-b-c \atop e',1+d-a-b, d+e-a-b};1\right],
\label{kdf4}
\eea
and holds under the extra condition $\Re(1-d+c+c')>0$.

So there remain (\ref{kdf1}) and (\ref{kdf2}) to prove. Our
proof of (\ref{kdf1}) is inspired by Gasper's proof~\cite{private} of the
$q$-analogue of one of our earlier results (eq.\ (40) of
Ref.~\cite{Vanderjeugt}), and uses the Beta function~:
\[
B(x,y)=\Gamma\left[ x,y \atop x+y \right]=\int_0^1
t^{x-1}(1-t)^{y-1}dt, \qquad\Re(x)>0, \Re(y)>0.
\]
Considering $B(x+m,y+n)$ and using $\Gamma(x+m)=(x)_m\Gamma(x)$,
one finds
\beq
{(x)_m(y)_n\over (x+y)_{m+n} } = \Gamma\left[ x+y\atop
x,y\right] \int_0^1 t^{x+m-1}(1-t)^{y+n-1}dt.
\label{beta}
\eeq
We start now from the series $F(1,1)$ with one constraint
$a'=d-a$ and apply~(\ref{beta}) to $(a)_m(d-a)_n/(d)_{m+n}$~:
\bea
F(1,1)&=&\sum_{m,n=0}^\infty {(a,b,c)_m(d-a,b',c')_n\over
(d)_{m+n} (e)_m (e')_n m!n!} \nn \\[1mm]
&=& \Gamma\left[d\atop a,d-a\right] \int_0^1 \sum_{m,n}
{(b,c)_m(b',c')_n \over (e)_m m! (e')_n n!}
t^{a+m-1}(1-t)^{d-a+n-1}dt. 
\label{tmp1}
\eea
In this last expression, we use Euler's identity (see
Ref.~\cite{Slater}, eq.\ (1.3.15))~:
\beq
\sum_m{(b,c)_m t^m\over (e)_m m!}={\ }_2F_1\left[{b,c \atop
e};t\right] =
(1-t)^{e-b-c}{\ }_2F_1\left[{e-b,e-c\atop e};t\right].
\eeq
We obtain
\beq
\Gamma\left[d\atop a,d-a\right] \int_0^1 \sum_{m,n}
{(e-b,e-c)_m(b',c')_n \over (e)_m m! (e')_n n!}
t^{a+m-1}(1-t)^{d+e-a-b-c+n-1}dt. 
\eeq
For the last integral, we apply again~(\ref{beta}). This leads
to 
\beq
\Gamma\left[d,d+e-a-b-c\atop d-a,d+e-b-c\right]
\sum_{m,n=0}^\infty { (a,e-b,e-c)_m (b',c',d+e-a-b-c)_n \over
(e)_m (e')_n (d+e-b-c)_{m+n} m!n!}.
\label{tmp2}
\eeq
So far, we have used only the condition $a'=d-a$. Suppose now a
second constraint is satisfied~: $e'=d+e-a-b-c$.
Then~(\ref{tmp2}) simplifies and using $(d+e-b-c)_{m+n} =
(a+e')_{m+n}= (a+e')_m(a+e'+m)_n$, we find that~(\ref{tmp2}) can
be written as
\beq
\Gamma\left[ d,e' \atop a',a+e'\right] \sum_{m=0}^\infty {
(a,e-b,e-c)_m \over (e,a+e')_m m!} {\ }_2F_1\left[{b',c'\atop
a+e'+m};1 \right].
\label{tmp3}
\eeq
Applying Gauss's theorem for the ${}_2F_1$ series yields, after
some elementary manipulations, the transformation
formula~(\ref{kdf1}). The extra condition
$\Re(a')>0$ comes from the convergence requirements of the
${}_4F_3$ series. The other conditions used to apply the Beta function
integral~(\ref{beta}) disappear by analytic continuation.

Next, consider again~(\ref{tmp2}) but now with the extra
constraint $d+e-a-b-c=1$. 
Then~(\ref{tmp2}) can be rewritten as
\beq
\Gamma\left[d\atop d-a\right]
\sum_{m=0}^\infty { (a,e-b,e-c)_m \over
(e)_m m!} {\ }\sum_{n=0}^\infty {(b',c')_n\over (e')_n \Gamma[1+a+m+n]}.
\label{tmp4}
\eeq
Consider now condition (i) where $a$ is a negative integer, say $a=-N$.
Then the
sum over $m$ is finite going from $0$ upto $N$, and due to the Gamma
function the summation over $n$ goes from $N-m$ upto $\infty$.
Replacing $n$ by $k+N-m$ and using $(x)_{k+N-m}=(x+N-m)_k \Gamma(x+N-m)
/ \Gamma(x)$ for $x=b',c',e'$, (\ref{tmp4}) reduces to
\bea
&&\Gamma\left[d,e'\atop d+N,b',c'\right]
\sum_{m=0}^{N} { (-N,e-b,e-c)_m \over (e)_m m!} \nn\\[1mm]
&&\qquad\times\Gamma\left[b'+N-m,c'+N-m\atop e'+N-m\right]
\sum_{k=0}^\infty {(b'+N-m,c'+N-m)_k\over (e'+N-m)_k k!}.
\label{tmp5}
\eea
The $k$-summation is a $_2F_1$ series and can be summed using Gauss's
theorem. Performing this explicitly, and using some elementary
manipulations with Gamma functions and Pochammer symbols, one arrives
at
\beq
\Gamma\left[ {d,e',b'+N,c'+N,-N-b'-c'+e' \atop
d+N,b',c',e'-b',e'-c'} \right]
{\ }_4F_3\left[{-N,e-b,e-c,-N-b'-c'+e' \atop 
e,1-N-b',1-N-c'};1\right].
\label{tmp6}
\eeq
The above ${}_4F_3$ series is finite because $N$ is a positive integer,
and thus one can apply reversal of series to it~:
\bea
&&{}_4F_3\left[{A,B,C,-N\atop D,E,F};1\right]=(-1)^N {(A,B,C)_N \over
(D,E,F)_N}\nn\\[1mm]
&&\times {\ }_4F_3\left[ {1-D-N,1-E-N,1-F-N,-N \atop 1-A-N,1-B-N,
1-C-N};1\right].
\eea
This gives rise to~(\ref{kdf2}),
thus proving it under the condition~(i).

In order to prove~(\ref{kdf2}) -- with $a'=d-a$ and $d+e-a-b-c=1$ -- 
under the conditions~(ii), i.e.\ $c=-N$ and
$c'=-N'$ with $N$ and $N'$ positive integers, we shall make use of the
following identity ($N$ positive integer), 
\beq
{}_3F_2\left[{A,B,-N \atop C,D};1\right] = {(C-A)_N \over (C)_N}
{\ }_3F_2\left[{A,D-B,-N \atop 1+A-C-N,D};1\right],
\label{Slatertf}
\eeq
which follows, for example, from eq.~(4.3.4.2) of Slater~\cite{Slater}.
Starting from~(\ref{kdf}) and using $(d)_{m+n}=(d)_n(d+n)_m$, $F(1,1)$
can be rewritten as
\beq
\sum_{n=0}^{N'} {\ }_3F_2\left[{a,b,-N\atop d+n,e};1\right]
{(a',b',-N')_n \over (d,e')_n n!}.
\label{tmp7}
\eeq
Applying (\ref{Slatertf}) to this $_3F_2$ leads to
\bea
{\ }_3F_2\left[{a,b,-N\atop d+n,e};1\right]&=&
{(d+n-a)_{N}\over (d+n)_{N}}{\ }_3F_2\left[{a,e-b,-N\atop
e-b-n,e};1\right]\nn \\
&=& {(d-a)_{N}\over (d)_{N}} {(d,d-a+N)_n\over
(d-a,d+N)_n}{\ }_3F_2\left[{a,e-b,-N\atop
e-b-n,e};1\right] .
\label{tmp8}
\eea
Plugging this in (\ref{tmp7}), it becomes
\beq
{(d-a)_{N}\over (d)_{N}}  \sum_{m,n}
{(a,-N,e-b)_m \over (e,e-b-n)_m m!} {(d-a+N,b',-N')_n \over (d+N,e')_n
n!}. 
\label{tmp9}
\eeq
Since $d+e-a-b-c=d+e-a-b+N=1$, we have
\beq
{(e-b)_m(d-a+N)_n\over (e-b-n)_m}=
{(e-b)_m(1+b-e)_n\over (e-b-n)_m}=(1-e+b-m)_n,
\eeq
and thus (\ref{tmp9}) reduces to
\beq
{(d-a)_{N}\over (d)_{N}}  \sum_{m}
{(a,-N)_m \over (e)_m m!}{\ }_3F_2\left[{b',1-e+b-m,-N'\atop e',d+N};1\right].
\label{tmp10}
\eeq
Since also $N'$ is a positive integer, we can apply~(\ref{Slatertf}) to
the last expression, and obtain
\beq
{(d-a)_{N}\over (d)_{N}}{(e'-b')_{N'}\over(e')_{N'}}  \sum_{m,n}
{(a,-N)_m \over (e)_m m!}{(b',a+m,-N')_n\over (1+b'-N'-e',d+N)_n n!}.
\label{tmp11}
\eeq
Making the replacement $(a)_m(a+m)_n=(a)_n(a+n)_m$ implies that the
last sum can be rewritten as
\beq
{(d-a)_{N}\over (d)_{N}}{(e'-b')_{N'}\over(e')_{N'}}  \sum_{n}
{(b',a,-N')_n\over (1+b'-N'-e',d+N)_n n!} {\ }_2F_1\left[{a+n,-N\atop
e};1 \right],
\label{tmp12}
\eeq
to which Vandermonde's theorem can be applied. Replacing 
$(e-a-n)_{N}$  by $(e-a)_{N}(1-e+a)_n/(1-e+a-N)_n$ leads to
\beq
{(d-a)_{N}\over (d)_{N}}{(e-a)_{N}\over(e)_{N}} 
{(e'-b')_{N'}\over(e')_{N'}}  \sum_{n}
{(b',a,-N',1-e+a)_n\over (1+b'-N'-e',d+N,1-e+a-N)_n n!} .
\label{tmp13}
\eeq
Using $d+e-a-b+N=1$ this can finally be rewritten in the form
(\ref{kdf2}), providing a proof under the condition (ii).

\section{Special cases and summation formulas}

Some limiting cases of (\ref{kdf1})--(\ref{kdf3}) are worth
considering. If we assume that there are three relations among
the 9 parameters in~(\ref{kdf}), i.e.
\beq
e'=d+e-a-b-c,\quad a'=d-a,\quad b'=d-b,
\eeq
then, using Dixon's theorem (eq.~(2.3.3.7) of
Ref.~\cite{Slater}) one can deduce from~(\ref{kdf1}) that
\bea
F(1,1)&=&\Gamma\left[ e,e' \atop e-c, e'+c \right]{\ }_3F_2\left[
{d-a, d-b, c+c' \atop d, e'+c};1\right],\nn\\[1mm]
&&\Re(e')>0,\qquad \Re(e-c-c')>0.
\label{3f2}
\eea
In particular, for $c'=-c$, the rhs simply reduces to a product
of Gamma functions, and we obtain the summation formula
\bea
&&F_{1:1}^{0:3}\left[ {\hy \atop d}:{a,b,c\atop
e};{d-a,d-b,-c\atop d+e-a-b-c};1,1\right]=
\Gamma\left[ e,e+d-a-b-c \atop e-c,e+d-a-b \right],
\label{res1} \\[1mm]
&&\qquad \Re(e)>0,\qquad \Re(d+e-a-b-c)>0.\nn
\eea
This special formula was already obtained
earlier~\cite{Vanderjeugt} in the context of 9-$j$
coefficients and proved by Karlsson~\cite{Karlsson2}; its
$q$-analogue was proved by Gasper~\cite{private}.

If on the other hand, we take in~(\ref{3f2}) the extra
condition $c'=d-c$, the ${}_3F_2$ series reduces to a ${}_2F_1$,
which can be summed with Gauss's theorem. There results the
following summation formula~:
\bea
&&F_{1:1}^{0:3}\left[ {\hy \atop d}:{a,b,c\atop
e};{d-a,d-b,d-c\atop d+e-a-b-c};1,1\right]=
\Gamma\left[ e,e+d-a-b-c,e-d \atop e-a,e-b,e-c \right],
\label{res2}\\[1mm]
&&\qquad \Re(e-d)>0,\qquad \Re(d+e-a-b-c)>0.\nn
\eea

Another interesting set of three relations among the 9
parameters is
\beq
e'=d+e-a-b-c,\quad a'=d-a,\quad d=b+c+b'+c'.
\eeq
The transformation formula (\ref{kdf1}) now becomes
\bea
F(1,1)&=&\Gamma\left[ d,e,e' \atop a',e+b', e+c' \right]{\ }_3F_2\left[
{a, e-b, e-c \atop e+b', e+c'};1\right],\nn\\[1mm]
&&\Re(e')>0,\qquad \Re(e)>0.
\label{3f22}
\eea
The summation formulas that can be deduced from here are
again~(\ref{res1}) and~(\ref{res2}). In principle one can
deduce further summation formulas by requiring the ${}_3F_2$
in~(\ref{3f2}) or~(\ref{3f22}) to be terminating and
Saalsch\"utzian or by requiring the above ${}_3F_2$'s or
${}_4F_3$'s to be of Karlsson-Minton
type~\cite{Karlsson1,Minton} (see also Eq.~(1.9.1)
of~\cite{Gasper}). We give one example~: consider~(\ref{3f2})
with $c'=e-c-1$ and $d-a$ or $d-b$ a negative integer. Then the
Pfaff-Saalsch\"utz formula can be used and one obtains~:
\bea
&&F_{1:1}^{0:3}\left[ {\hy \atop d}:{a,b,c\atop
e};{d-a,d-b,e-c-1\atop d+e-a-b-c};1,1\right]= \nn\\[1mm]
&&\qquad \Gamma\left[ 1-a,1-b,e,e-d,d+e-a-b-c \atop
1-d,e-a,e-b,e-c,1+d-a-b \right], \nn\\[1mm]
&&\Re(d+e-a-b-c)>0\quad\hbox{and}\quad d-a\ \hbox{or}\ d-b\ 
\hbox{a negative integer}. \label{res3}
\eea

Some special summation formulas follow by specializing (\ref{kdf2})
under one of the two conditions. In the case (i), one can specialize
$b'=d-b$ and $c'=d-c$; or $b'=d-b$ and $e'=1+c+c'$. In the case (ii),
one can choose $b'=d-b$ and $e'=1-a-b+c'+d$; or $b'=d-b$ and
$e'=1+c+c'$. In terms of 5 independent parameters, such specializations
give rise to the following four summation formulas~:
\bea
&&F_{1:1}^{0:3}\left[ {\hy \atop d}:{-N,b,c\atop
1-N+b+c-d};{d+N,d-b,d-c\atop e'};1,1\right]= \nn\\[1mm]
&&\qquad {(d-b,d-c,1+d-e')_{N}\over(d,d-b-c,1+2d-b-c-e')_{N}}
\Gamma\left[ e',e'+b+c-2d \atop e'+b-d,e'+c-d \right], \nn\\[1mm]
&&\Re(e'-N+b+c-2d)>0\quad\hbox{and}\quad N\ \ 
\hbox{a positive integer}; \label{g1}
\eea
\bea
&&F_{1:1}^{0:3}\left[ {\hy \atop d}:{-N,b,c\atop
1-N+b+c-d};{d+N,d-b,c'\atop 1+c+c'};1,1\right]= \nn\\[1mm]
&&\qquad {(d-b,d-c-c')_{N}\over(d,d-b-c)_{N}}
\Gamma\left[ 1+c+c',1+b+c-d \atop 1+c,1+b+c+c'-d \right], \nn\\[1mm]
&&\Re(1-N+b+c-d)>0\quad\hbox{and}\quad N\ \ 
\hbox{a positive integer}; \label{g2}
\eea
\bea
&&F_{1:1}^{0:3}\left[ {\hy \atop d}:{a,b,-N\atop
1+a+b-N-d};{d-a,d-b,-N'\atop 1-a-b-N'+d};1,1\right]= \nn\\[1mm]
&&\qquad {(d-a,d-b)_{N}(a,b)_{N'}\over(d)_{N+N'}(d-a-b)_{N}
(a+b-d)_{N'} }, \nn\\[1mm]
&&N\quad\hbox{and}\quad N'\ \ 
\hbox{positive integers}; \label{fi1}
\eea
\bea
&&F_{1:1}^{0:3}\left[ {\hy \atop d}:{a,b,-N\atop
1+a+b-N-d};{d-a,d-b,-N'\atop 1-N-N'};1,1\right]= \nn\\[1mm]
&&\qquad {(d-a,d-b)_{N+N'}\over(d)_{N+N'}(d-a-b)_{N}
(N)_{N'} }, \nn\\[1mm]
&&N(\ne 0)\quad\hbox{and}\quad N'\ \ 
\hbox{positive integers}. \label{fi2}
\eea

\section*{Acknowledgements}

It is a pleasure to thank Professor G.\ Gasper for some
fruitful exchanges, and particularly for pointing out a major error in
the first version of the manuscript.
This research was partly supported by the E.E.C. (contract No.
CI1*-CT92-0101).


\end{document}